\documentclass[12pt]{amsart}
\usepackage[colorlinks=true,citecolor=blue,linkcolor=blue]{hyperref}
\usepackage{latexsym,amsmath,amssymb}
\usepackage{cleveref}
\usepackage{accents}

\title{\protect{On BMO and Carleson measures on Riemannian Manifolds}}

plus 6pt minus 12pt
\def\M{{\mathcal M}}
\def\N{{\mathbb N}}

\newtheorem{theorem}{Theorem}
\newtheorem{lemma}[theorem]{Lemma}
\newtheorem{corollary}[theorem]{Corollary}
\newtheorem{proposition}[theorem]{Proposition}

\newtheorem{remark}[theorem]{Remark}
\newtheorem{definition}[theorem]{Definition}

\def\supp{{\rm supp\,}}

\newcommand{\R}{\mathbb{R}}

\newcommand{\brac}[1]{\left (#1 \right )}

\author[D. Brazke] {Denis Brazke}
\address[D. Brazke]{Department of Mathematics,
University of Heidelberg,
Im Neuenheimer Feld 205,
69120 Heidelberg, Germany}
\email{brazke@stud.uni-heidelberg.de}

\author[A. Schikorra]{Armin Schikorra}
\address[A. Schikorra]{Department of Mathematics,
University of Pittsburgh,
301 Thackeray Hall,
Pittsburgh, PA 15260, USA}
\email{armin@pitt.edu}

\author[Y. Sire]{Yannick Sire}
\address[Y. Sire]{Johns Hopkins University, Krieger Hall, Baltimore, USA}
\email{sire@math.jhu.edu}

\newcommand{\barint}{
\rule[.036in]{.12in}{.009in}\kern-.16in \displaystyle\int }

\newcommand{\barcal}{\mbox{$ \rule[.036in]{.11in}{.007in}\kern-.128in\int $}}

\def\mvint_#1{\mathchoice
          {\mathop{\vrule width 6pt height 3 pt depth -2.5pt
                  \kern -8pt \intop}\nolimits_{\kern -3pt #1}}%
          {\mathop{\vrule width 5pt height 3 pt depth -2.6pt
                  \kern -6pt \intop}\nolimits_{#1}}%
          {\mathop{\vrule width 5pt height 3 pt depth -2.6pt
                  \kern -6pt \intop}\nolimits_{#1}}%
          {\mathop{\vrule width 5pt height 3 pt depth -2.6pt
                  \kern -6pt \intop}\nolimits_{#1}}}

\numberwithin{theorem}{section} \numberwithin{equation}{section}

\newcommand{\lap}{\Delta }
\newcommand{\aleq}{\lesssim}

\newcommand{\aeq}{\approx}

\newcommand{\laps}[1]{(-\lap)^{\frac{#1}{2}}}

\begin{document}
\begin{abstract}
Let $\mathcal{M}$ be a Riemannian $n$-manifold with a metric such that the manifold is Ahlfors-regular. We also assume either non-negative Ricci curvature, or that the Ricci curvature is bounded from below together with a bound on the gradient of the heat kernel. We characterize BMO-functions $u: \mathcal{M} \to \R$ by a Carleson measure condition of their $\sigma$-harmonic extension $U: \mathcal{M} \times (0,\infty) \to \R$. We make crucial use of a $T(b)$ theorem proved by Hofmann, Mitrea, Mitrea, and Morris.

As an application we show that the famous theorem of Coifman--Lions--Meyer--Semmes holds in this class of manifolds: Jacobians of $W^{1,n}$-maps from $\mathcal{M}$ to $\R^n$ can be estimated against BMO-functions, which now follows from the arguments for commutators  recently proposed by Lenzmann and the second-named author using only harmonic extensions, integration by parts, and trace space characterizations.
\end{abstract}

\maketitle

\tableofcontents

\section{Introduction}
It is a classical result that in Euclidean space there is a relation between $BMO$-functions $u: \R^n\to \R$ and Carleson measures in $\R^{n+1}_+$. Precisely, the following statement can be found, e.g., in \cite[IV, \textsection 4.3, Theorem~3, pp.159 and \textsection 4.4.3]{Stein-1993}.
\begin{theorem}\label{th:carlesonclassic}
Let $u \in C_c^\infty(\R^n)$ and denote by $U(x,t): \R^{n+1}_+ \to \R$ the harmonic extension, i.e. the unique solution to the following equation:
\begin{equation}\label{eq:he}
\begin{cases}
\lap_{x,t} U \equiv (\partial_{tt} + \lap_x) U = 0 \quad &\mbox{in $\R^{n+1}_+$}\\
U = u &\mbox{on $\R^n \times \{0\}$}\\
\lim_{|(x,t)| \to \infty} U(x,t) = 0
\end{cases}
\end{equation}
Then the following two $BMO$-seminorms are equivalent: The integral one
\[
 [u]_{BMO} = \sup_{B \subset \R^n} |B|^{-1} \int_{B}|u-(u)_{B}|,
\]
and the Carleson-measure version
\[
 [u]_{\tilde{BMO}} := \brac{\sup_{B \subset \R^n} |B|^{-1} \int_{T(B)} t|\nabla_{(x,t)} U(x,t)|^2\, \mathrm dx\,\mathrm dt }^{\frac{1}{2}}.
\]
Here, $(u)_B = |B|^{-1} \int_B u$ and $T(B)$ is the tent in $\R^n \times (0,\infty)$ over the ball $B$, namely if $B = B(x_0,r)$ then $T(B) = \{(x,t): |x-x_0| < r-t\}$.
\end{theorem}

While relations between Carleson measures and \emph{certain} extensions of functions have been extended to spaces of homogeneous type (see e.g. \cite{HS08,TT16,DMC17}),  these extensions are usually of a potential type (with conditions on kernel decay). The main drawback is that these extensions in general  do not satisfy an equation such as \eqref{eq:he}. On the other hand, 
in applications such as proving sharp commutator estimates, see \cite{LS}, it is beneficial (and maybe even crucial) to have the extension satisfying certain PDEs such as \eqref{eq:he} (or more generally satisfying a Dirichlet-to-Neumann principle \cite{CS07}). 

The aim of the present work is to prove an equivalence result like in the previous theorem involving a natural PDE-extension to the half-space in a rather general geometric framework.  Let $(\mathcal{M},g)$ be an $n$-dimensional smooth Riemannian manifold which is also Ahlfors regular, meaning that the measure of a ball of radius $r$ is (uniformly) comparable to $r^n$. By $\lap_{\mathcal{M}}$ we denote the Laplace-Beltrami operator on $\mathcal{M}$. We equip $\M$ with the Carnot-Carath\'eodory metric $d$. Without further assumptions on the manifold it seems implausible that, e.g., harmonic extensions satisfy a statement such as \Cref{th:carlesonclassic}. In Theorem~\ref{main} we introduce such assumptions on the manifold and its heat kernel.

Besides the classical harmonic extension we also take $\sigma$-harmonic extensions into account. To define them we follow the semigroup representation (cf. \cite{Stinga-Torrea-2010} for the Euclidean analogue but as stated in the latter it extends to much more general contexts).

Let us clarify the differential operators that we use. Denote by $\mathrm d$  the exterior derivative and $\star$ the Hodge operator. We define the Laplace-Beltrami operator (or Laplace-de-Rham operator, which are the same for us since they act only on functions/0-forms) by $\Delta = \Delta_\M := \star  \, \mathrm d \star \mathrm d$, and the gradient of a smooth function $f$ by $\nabla f = \nabla_\M f= \mathrm d f$, (or depending on the context $\nabla f = (\mathrm d f)^\sharp$). With this setup, we have (with standard abuse of notation) $|\nabla f|^2 = g(\nabla f,\nabla f)$ and
\[
		\langle \nabla f,\nabla h \rangle_{L^2} = \int_\M \star \mathrm d f \wedge \mathrm d h
		\]
		\[
		\text{i.e. } \|\nabla f\|_{L^2}^2 = \int_\M g(\nabla f,\nabla f) \ \mathrm dx = \int_\M |\nabla f|^2 \ \mathrm dx.
		\]
The $\sigma$-harmonic extension is defined as follows.
\begin{definition}[$\sigma$-Harmonic extension to $\mathcal{M} \times (0,\infty)$]\label{defExtension} Let $\mathcal{M}$ be as above and $0 < \sigma < 1$.
For $u \in C_c^\infty(\mathcal{M})$ the $\sigma$-harmonic extension $U: \mathcal{M} \times [0,\infty) \to \R$ is the solution to                                                                                                                                                          
\[
 \begin{cases}
 \lap_{\mathcal{M}} U  + \frac{1-2\sigma}{t}\partial_t U + \partial_{tt} U = 0 \quad& \mbox{in $\mathcal{M}\times (0,\infty)$}\\
 U(x,0) = u(x) & \mbox{in $\mathcal{M}$}\\
\lim_{|(x,t)| \to \infty} U(x,t) = 0. 
 \end{cases}
\]
This solution is formally given by
\[
 U(x,t) = \frac{1}{4^\sigma \Gamma(\sigma)} t^{2\sigma} \int_0^\infty e^{s\lap_{\mathcal{M}}} u(x)\, e^{-\frac{t^2}{4s}}\, \frac{\mathrm ds}{s^{1+\sigma}},
\]
and explicitly one has
\[
\begin{split}
 U(x,t) =& \frac{1}{4^\sigma \Gamma(\sigma)} \int_0^\infty \int_{\mathcal{M}} p(x,y,s) u(y) \mathrm dy\, t^{2\sigma} e^{-\frac{t^2}{4s}}\, \frac{\mathrm ds}{s^{1+\sigma}},	\\
 			= & \frac{1}{\Gamma(\sigma)} \int_0^\infty \int_\M p(x,y,\tfrac{t^2}{4s}) \, u(y) \, \mathrm dy \, e^{-s} \, s^{\sigma - 1} \, \mathrm ds,
\end{split}
 \]
where $p(x,y,s)$ is the heat kernel for $\mathcal{M}$, i.e. 
\begin{equation}\label{eq:heatkerneleq}
\begin{cases} 
(\partial_t - \lap_x)p(x,y,s) = 0 \quad &\mbox{for all $x,y \in \mathcal{M}$ and $s > 0$}\\
            p(x,y,0) = \delta_{x,y}                                                                                                                            . \end{cases}
                                                                                                                                                                \end{equation}
   
\end{definition}
The previous definition is not explicitly stated in \cite{Stinga-Torrea-2010} but it is easy to check that the semi-group approach automatically carries over to such a geometric setting under very weak assumptions on the manifold, see Section~\ref{a:heatkernelextension}. For more information and properties about the heat kernel, see \cite{Grig03}.
 	
We define the following semi-norms: Let $U(x,t)$ be the $\sigma$-harmonic extension of $u$ to $\mathcal{M} \times (0,\infty)$. Denote the usual $BMO$-norm as
\[
 [u]_{BMO(\mathcal{M})} := \sup_{B \subset \mathcal{M}} |B|^{-1} \int_{B}|u-(u)_{B}|,
\]
where the supremum is taken over balls $B$. Furthermore, we define a notion of $BMO$ in terms of the $\sigma$-harmonic extension and Carleson measures, namely
\begin{equation}\label{eq:carlesonbmo}
 [u]_{\tilde{BMO}(\mathcal{M})} \equiv [u]_{\tilde{BMO}_\sigma(\mathcal{M})}:= \brac{\sup_{B \subset \mathcal{M}} |B|^{-1} \int_{T(B)} t|\nabla_{x,t} U(x,t)|^2\, \mathrm dx\, \mathrm dt }^{\frac{1}{2}}.
\end{equation}
Again, $T(B)$ is the tent in $\mathcal{M} \times (0,\infty)$ over the ball $B$, namely if $B = B(x_0,r)$ then $T(B) = \{(x,t): d(x,x_0) < r-t\}$.

Our main result is the following:
\begin{theorem}\label{main}
Let $\M$ be a complete path-connected and Ahlfors regular manifold without boundary, such that the Ricci curvature is bounded from below. 

%
If moreover the heatkernel of $\mathcal{M}$ satisfies 
\[
 \sup_{x,y \in \mathcal{M}} |\nabla p(x,y,t)| \aleq t^{-\frac{n+1}{2}},
\]
then for any $0 < \sigma < 1$ the semi-norms of $BMO$ defined above are equivalent, i.e  for any $u \in C_c^\infty(\mathcal{M})$ we have 
\[
 [u]_{BMO(\mathcal{M})}  \aeq [u]_{\tilde{BMO}(\mathcal{M})}.
\]
\end{theorem}

In case of non-negative Ricci curvature of the manifold, we can drop the assumption of the gradient bound and have the following:

\begin{theorem}\label{main2}
Let $\M$ be a complete path-connected and Ahlfors regular manifold without boundary, such that the Ricci curvature is non-negative. Then for any $0 < \sigma < 1$ the semi-norms of $BMO$ defined above are equivalent, i.e  for any $u \in C_c^\infty(\M)$ we have 
\[
[u]_{BMO(\M)}  \approx [u]_{\tilde{BMO}(\M)}.
\]
\end{theorem}

Theorem~\ref{main} and Theorem~\ref{main2} can be very useful for example when estimating commutators via harmonic extensions, as recently proposed in \cite{LS}, who gave new proofs for a large class of commutator estimates (in Euclidean space). Their argument is based on integration by parts and trace space characterizations for $\sigma$-harmonic extensions. Since in this paper we obtained the latter characterization for BMO, one can follow the ideas in \cite{LS} almost verbatim for manifolds. For example, the following estimate on Jacobians was obtained for $\mathcal{M}=  \R^n$ in the celebrated work \cite{CLMS}, which lead to several breakthroughs in regularity theory e.g. of harmonic maps. We can extend it to manifolds.
\begin{theorem}\label{th:jac}
Let $\mathcal{M}$ be an $n$-manifold as in Theorem~\ref{main} or Theorem~\ref{main2}, $n \geq 2$.
For $f \in W^{1,n}(\mathcal{M},\R^n)$ and any $\varphi \in C_c^\infty(\mathcal{M})$ we have
\[
 \int_{\mathcal{M}} \mathrm df^1\wedge \ldots \wedge \mathrm df^n\ \varphi \leq C(\mathcal{M}) \ \|\nabla f\|_{L^n(\M)}^n\ [\varphi]_{BMO}.
\]
\end{theorem}
The remainder of this paper will be as follows: In Section~\ref{s:admissible} we introduce the notion of an admissible manifold, which is more general than the one in Theorem~\ref{main} or Theorem~\ref{main2}, but more complicated to check. In Section~\ref{s:Tb} we use the $T(b)$-theorem from \cite{Hofmann-Mitreas2017} to obtain square function estimates. In Section~\ref{s:mainproof} we prove Theorem~\ref{main} and Theorem~\ref{main2}. In Section~\ref{s:thjac} we prove the Jacobian estimate, Theorem~\ref{th:jac}. Computations concerning the $\sigma$-harmonic extensions are moved to the appendix, Section~\ref{a:heatkernelextension}.

{\bf Acknowledgments} 
The authors would like to warmly thank Laurent Saloff-Coste for highly valuable discussions on the topic of this paper and Fabrice Baudoin and Ryan Alvarado for providing crucial references. Part of this work was carried out while D.B. was visiting the University of Pittsburgh, he likes to thank the Math Department for its hospitality.

Y.S. acknowledges funding from Simons Foundation. A.S. acknowledges funding from the Simons foundation and the Daimler and Benz foundation.

\section{Admissible manifolds and heat kernel estimates}\label{s:admissible}
The proofs of \Cref{main} and \Cref{main2} are based on heat kernel estimates which allow to deduce a $T(b)$ theorem, see Section~\ref{s:Tb}. We exhibit a large class of Riemannian manifolds for which our theorem applies.  We first define a general setting in which our theory actually works. 
\begin{definition}[Admissible Manifolds]\label{mfdadmissible}
A manifold $\M$ is said to be admissible, if it is complete, path-connected, Ahlfors regular, without boundary and its heat kernel $p(x,y,t)$ satisfies the following conditions:
For every $\sigma > 0$ there exists $\nu > 0$ such that for all $t > 0$, $x, y \in \mathcal{M}$:
\begin{equation}\label{c1}
 \int_0^\infty p(x,y,t^2 s) \brac{1+\frac{1}{s}} e^{-\frac{1}{4s}}\, \frac{\mathrm ds}{s^{1+\sigma}} \aleq \frac{t^{\nu}}{\brac{d(x,y)^2+t^2}^{\frac{n+\nu}{2}}},
\end{equation}
\begin{equation}\label{c2}
 \int_0^\infty |\nabla_xp(x,y,t^2 s)| \, e^{-\frac{1}{4s}}\, \frac{\mathrm ds}{ s^{1+\sigma}}\aleq \frac{t^{\nu-1}}{\brac{d(x,y)^2+t^2}^{\frac{n+\nu}{2}}},
\end{equation}
\begin{equation}\label{c3}
 \int_0^\infty |\nabla_y p(x,y,t^2 s)| \brac{1+\frac{1}{s}} e^{-\frac{1}{4s}}\, \frac{\mathrm ds}{s^{1+\sigma}}\aleq \frac{t^{\nu-1}}{\brac{d(x,y)^2+t^2}^{\frac{n+\nu}{2}}},
\end{equation}
\begin{equation}\label{c4}
 \int_0^\infty |\nabla_x \nabla_y p(x,y,t^2 s)| \, e^{-\frac{1}{4s}}\, \frac{\mathrm ds}{ s^{1+\sigma}} \aleq \frac{t^{\nu-2}}{\brac{d(x,y)^2+t^2}^{\frac{n+\nu}{2}}}.
\end{equation}
\end{definition}

The following lemma is providing a more treatable class of admissible manifolds. 
\begin{lemma}\label{la:admissiblenicer}
Let $\M$ be a complete path-connected and Ahlfors regular manifold without boundary. 
If the heat kernel satisfies
	\begin{equation}\label{eq:heat1}
				p(x,y,t) \lesssim t^{-\frac n2} \, e^{-c \frac{d(x,y)^2}{t}},
			\end{equation}
			\begin{equation}\label{eq:heat2}
				|\nabla p(x,y,t)| \lesssim t^{-\frac{n+1}{2}}\, e^{-c \frac{d(x,y)^2}{t}},
			\end{equation}
			\begin{align} \label{eq:heat3}
				|\nabla_{\kern -2pt x}\nabla_{\kern -2pt y} p(x,y,t)| \lesssim \frac{1}{t^{\frac n2 + 1}} e^{-c\frac{d(x,y)^2}{t}},
			\end{align}
for some constant $c > 0$, then $\M$ is admissible with $\nu = 2\sigma$.
\end{lemma}
\begin{proof} The geometry of the manifold is the same as in Definition \ref{mfdadmissible}, so we only have to check the conditions \eqref{c1} to \eqref{c4}.

\underline{As for \eqref{c1}:}
Using the change of variables, we see that
		\begin{align*}
			\int_0^\infty p(x,y,t^2 s) \brac{1+\frac{1}{s}} e^{-\frac{1}{4s}}\, \frac{\mathrm ds}{s^{1+\sigma}}	& = t^{2 + 2\sigma - 2}\int_0^\infty p(x,y,s) \brac{1+\frac{t^2}{s}} e^{-\frac{t^2}{4s}}\, \frac{\mathrm ds}{s^{1+\sigma}}	\\
																											& \lesssim t^{2\sigma} \int_0^\infty \brac{1+\frac{t^2}{s}} e^{-c \frac{d(x,y)^2}{s}} e^{-c\frac{t^2}{s}} \frac{\mathrm ds}{s^{\frac n2 + 1 + \sigma}}.
		\end{align*}
	In the last step we employed \eqref{eq:heat1}. Set $\Lambda := d(x,y)^2 + t^2$, and making the change of variables $s \longmapsto \Lambda s$, we see
		\begin{align*}
			\ldots	& = t^{2 \sigma} \int_0^\infty \brac{1+\frac{t^2}{\Lambda s}} e^{-\frac cs} \frac{\mathrm ds}{s^{\frac n2 + 1 + \sigma}} \, \Lambda^{\frac n2 + \sigma} \lesssim \frac{t^{2\sigma}}{\Lambda^{\frac{n + 2\sigma}{2}}}  \int_0^\infty \brac{1+\frac{1}{s}} e^{-\frac cs} \frac{\mathrm ds}{s^{\frac n2 + 1 + \sigma}}.
		\end{align*}
	The integral can be estimated by a multiple of $\Gamma(\frac n2 + \sigma) + \Gamma(\frac n2 + \sigma + 1)$, so it is finite. This shows \eqref{c1}.

\underline{As for \eqref{c2} and for \eqref{c3}:} It suffices to show \eqref{c3}, since in our setting $\nabla_{\kern-2pt x} p(x,y,t) = \nabla_{\kern-2pt y} p(x,y,t)$ by the symmetry of the heat kernel. Using the gradient estimate \eqref{eq:heat2} we deduce
		\begin{align*}
			\int_0^\infty |\nabla_{\kern-2pt y} p(x,y,t^2 s)| \brac{1+\frac{1}{s}} e^{-\frac{1}{4s}}\, \frac{\mathrm ds}{s^{1+\sigma}}	& = t^{2 + 2\sigma - 2}\int_0^\infty |\nabla_{\kern-2pt y} p(x,y,s)| \brac{1+\frac{t^2}{s}} e^{-\frac{t^2}{4s}}\, \frac{\mathrm ds}{s^{1+\sigma}}	\\
																												& \lesssim t^{2\sigma} \int_0^\infty \brac{1+\frac{t^2}{s}} e^{-c \frac{d(x,y)^2}{s}} e^{-c\frac{t^2}{s}} \frac{\mathrm ds}{s^{\frac {n + 1}2 + 1 + \sigma}}.	
		\end{align*}
	As in \eqref{c1}, setting $\Lambda := d(x,y)^2 + t^2$, and making the change of variables $s \longmapsto \Lambda s$, we see
		\begin{align*}
			\ldots	& = t^{2 \sigma} \int_0^\infty \brac{1+\frac{t^2}{\Lambda s}} e^{-\frac cs} \frac{\mathrm ds}{s^{\frac {n + 1}2 + 1 + \sigma}} \, \Lambda^{\frac {n + 1}2 + \sigma} \lesssim \frac{t^{2\sigma}}{\Lambda^{\frac{n + 1 + 2\sigma}{2}}}  \int_0^\infty \brac{1+\frac{1}{s}} e^{-\frac cs} \frac{\mathrm ds}{s^{\frac {n + 1}2 + 1 + \sigma}}.
		\end{align*}
	The integral can be estimated by a multiple of $\Gamma(\frac {n + 1}2 + \sigma) + \Gamma(\frac {n + 1}2 + \sigma + 1)$, so it is finite. Moreover, since $\Lambda \geq t^2$, we get in the end
		\begin{align*}
			\int_0^\infty |\nabla_{\kern-2pt y} p(x,y,t^2 s)| \brac{1+\frac{1}{s}} e^{-\frac{1}{4s}}\, \frac{\mathrm ds}{s^{1+\sigma}} \lesssim  \frac{t^{2\sigma - 1}}{\Lambda^{\frac{n + 2\sigma}{2}}}.
		\end{align*}
	 This shows \eqref{c3}, and hence also \eqref{c2}. 

\underline{As for \eqref{c4}:} We proceed as in \eqref{c1} to \eqref{c3}. Using the estimate \eqref{eq:heat3}, we get
		\begin{align*}
			\int_0^\infty |\nabla_{\kern-2pt x} \nabla_{\kern-2pt y} p(x,y,t^2 s)| \, e^{-\frac{1}{4s}}\, \frac{\mathrm ds}{s^{1+\sigma}}	& = t^{2 + 2\sigma - 2}\int_0^\infty |\nabla_{\kern-2pt x} \nabla_{\kern-2pt y} p(x,y,s)| \, e^{-\frac{t^2}{4s}}\, \frac{\mathrm ds}{s^{1+\sigma}}	\\
																																	& \lesssim t^{2\sigma} \int_0^\infty  e^{-c \frac{d(x,y)^2}{s}} e^{-c\frac{t^2}{s}} \frac{\mathrm ds}{s^{\frac {n + 2}2 + 1 + \sigma}}.	
		\end{align*}
	As before, setting $\Lambda := d(x,y)^2 + t^2$, and making the change of variables $s \longmapsto \Lambda s$, we see
		\begin{align*}
			\ldots	& = t^{2 \sigma} \int_0^\infty e^{-\frac cs} \frac{\mathrm ds}{s^{\frac {n + 2}2 + 1 + \sigma}} \, \Lambda^{\frac {n + 2}2 + \sigma}	\lesssim \frac{t^{2\sigma}}{\Lambda^{\frac{n + 2 + 2\sigma}{2}}}  \int_0^\infty  e^{-\frac cs} \frac{\mathrm ds}{s^{\frac {n + 2}2 + 1 + \sigma}}.
		\end{align*}
	The integral can be rewritten as a multiple of $\Gamma(\frac {n + 2}2 + \sigma)$, so it is finite. Moreover, since $\Lambda \geq t^2$, we get in the end
		\begin{align*}
			\int_0^\infty |\nabla_{\kern-2pt x} \nabla_{\kern-2pt y} p(x,y,t^2 s)| e^{-\frac{1}{4s}}\, \frac{\mathrm ds}{s^{1+\sigma}} \lesssim  \frac{t^{2\sigma - 2}}{\Lambda^{\frac{n + 2\sigma}{2}}}.
		\end{align*}
	This finishes the proof.
\end{proof}

\begin{corollary}\label{co:2}
Let $\M$ be as in \Cref{main}. Then $\M$ is admissible.
\end{corollary}

\begin{proof}
The statement follows from known zero-order and first-order bounds on the heat kernel on $\M$ which we recall below.

By the curvature assumption we have \eqref{eq:heat1}, see \cite[Corollary 3.1]{LiYau}, see also \cite[Theorem 2.34]{baudoin}.
				
By \cite[Theorem 4.9]{coulhonSikora} the assumption on the heat kernel together with \eqref{eq:heat1} implies
			\[
				|\nabla p(x,y,t)| \lesssim t^{-\frac{n+1}{2}} \brac{1+\frac{d^2(x,y)}{t}}\, e^{-c \frac{d(x,y)^2}{t}}.
			\]
Since $(1+|x|) e^{-|x|} \leq C e^{-\frac{1}{2}|x|}$ this readily implies \eqref{eq:heat2}. Recall that the heat kernel is symmetric, so the gradient estimate holds both for $\nabla_x$ and $\nabla_y$.

For \eqref{eq:heat3} we use the semi-group property of the heat kernel, i.e.
			\begin{align*}
				p(x,y,2t) = \int_\M p(x,z,t) \, p(y,z,t) \, \mathrm dz.
			\end{align*}
		Keep in mind, that if $a \leq b + c$, then $a^2 \leq 2(b^2 + c^2)$. 
		So using H\"older's inequality, we arrive at
			\begin{align} \label{eq:saloff-coste} \begin{aligned}
				e^{\frac{d(x,y)^2}{At}} \, |\nabla_{\kern -2pt x}\nabla_{\kern -2pt y} p(x,y,2t)|	& \leq \int_\M e^{2\frac{d(x,z)^2}{At}} \, |\nabla_{\kern -2pt x} p(x,z,t)| \, e^{2\frac{d(y,z)^2}{At}} \, |\nabla_{\kern -2pt y} p(y,z,t)| \, \mathrm dz	\\
																									& \leq \Big( \int_\M e^{4\frac{d(x,z)^2}{At}} \, |\nabla_{\kern -2pt x} p(x,z,t)|^2 \, \mathrm dz \Big)^\frac 12 \Big( \int_\M e^{4\frac{d(y,z)^2}{At}} \, |\nabla_{\kern -2pt y} p(y,z,t)|^2 \, \mathrm dz \Big)^\frac 12	\\
																									& \lesssim \Big( \int_\M e^{4\frac{d(x,z)^2}{At}} \, \frac{1}{t^{n + 1}} \, e^{-c\frac{d(x,z)^2}{t}} \, \mathrm dz \Big)^\frac 12	\\
																									& \quad \times \Big( \int_\M e^{4\frac{d(y,z)^2}{At}} \, \frac{1}{t^{n + 1}} \, e^{-c\frac{d(y,z)^2}{t}} \, \mathrm dz \Big)^\frac 12	\\
																									& = \frac{1}{t^{n + 1}} \int_\M e^{4\frac{d(x,z)^2}{At}} \, e^{-c\frac{d(x,z)^2}{t}} \, \mathrm dz.
			\end{aligned} \end{align}
		We now choose $A$ so big, such that $\frac 4A < c$ in the equation above, which means we look for an estimate of the form
			\begin{align*}
				\int_\M e^{-c\frac{d(x,z)^2}{t}} \, \mathrm dz \lesssim t^{\frac n2},
			\end{align*}
		where $c > 0$. Let $B_0 = B(x,\sqrt t)$ and let $B_k = B(x, 2^k \sqrt t) \setminus B(x,2^{k - 1} \sqrt t)$. Then we have $\M = \cup B_k$. Furthermore, it holds
			\begin{align*}
				\int_{B_0} e^{-c\frac{d(x,z)^2}{t}} \, \mathrm dz \lesssim |B_0| \lesssim t^\frac n2
			\end{align*}
		by the Ahlfors regularity. Moreover, it holds
			\begin{align*}
				\int_{B_k} e^{-c\frac{d(x,z)^2}{t}} \, \mathrm dz \lesssim 2^{nk} \, e^{-c \, 2^{k - 1}} \, t^{\frac n2}
			\end{align*}
		again by the Ahlfors regularity. Together, we have
			\begin{align*}
				\int_\M e^{-c\frac{d(x,z)^2}{t}} \, \mathrm dz \leq \sum_{k = 0}^\infty \int_{B_k} e^{-c\frac{d(x,z)^2}{t}} \, \mathrm dz \lesssim t^{\frac n2} \sum_{k = 0}^\infty 2^{nk} \, e^{-c \, 2^{k - 1}} \lesssim t^{\frac n2}.
			\end{align*}
		This estimate together with \eqref{eq:saloff-coste} gives \eqref{eq:heat3} as desired.
\end{proof}

\begin{corollary}\label{co:3}
Let $\M$ be as inTheorem~\ref{main2}. Then $\M$ is admissible.
\end{corollary}

\begin{proof}
Using \cite[Theorem 4.2]{baud13}, we obtain the gradient estimate \eqref{eq:heat2}. The claim then follows as in the proof of Corollary~\ref{co:2}.
\end{proof}

\begin{remark}
The previous results, thanks to \cite{ACDH}, extend straightforwardly to Lie groups with polynomial volume (notice that these are spaces of homogeneous type in the sense of Coifman and Weiss \cite{CW77}). 
\end{remark}
\section{\texorpdfstring{$T(b)$}{T(b)}-Theorem and square function estimates}\label{s:Tb}
We use the following important version of the local $T(b)$-theorem on manifolds, which is proven in much greater generality in \cite[Theorem~3.7.]{Hofmann-Mitreas2017}. It allows us to pass from local estimates in small balls of $\mathcal{M}$ (i.e. essentially the Euclidean space) to global estimates.
\begin{theorem}\label{Tb}%
	
Let $\M$ be an admissible manifold and let $T$ be an operator, acting on functions $f:\M \to \R$ via
\[
 Tf(x,t) := \int_{\mathcal{M}} \kappa(x,y,t) \, f(y) \, \mathrm dy,
\]
where $\kappa : \M \times \M \times (0,\infty) \longrightarrow \R$ is integrable and satisfies 
\[
 \int_{\M} \kappa(x,y,t) \, \mathrm dy = 0 \quad \text{for all } x \in \M, \ t > 0,
\]
\begin{equation}\label{eq:TBkappaest}
 |\kappa(x,y,t)| \aleq \frac{t^{\nu}}{\brac{d(x,y)^2+t^2}^{\frac{n+\nu}{2}}} \quad  \text{for all } x,y \in \M, \ t > 0
\end{equation}
\begin{equation}\label{eq:TBkappalipest}
 \frac{|\kappa(x,y_1,t)-\kappa(x,y_2,t)|}{d(y_1,y_2)} \aleq \frac{t^{\nu - 1}}{\brac{d(x,y_1)^2+t^2}^{\frac{n+\nu}{2}}} \quad \text{for all } d(y_1,y_2) \leq \frac 12 (d(x,y_1)^2 + t^2)^\frac 12.
\end{equation}
Then,
\[
 \brac{\int_{\mathcal{M}} \int_0^\infty |Tf(x,t)|^2\, \frac{\mathrm dt}{t} \mathrm dx}^{\frac{1}{2}} \aleq \brac{\int_{\mathcal{M}} |f(x)|^2\, \mathrm dx}^{\frac{1}{2}}.
\]
\end{theorem}
\begin{proof}
All conditions in \cite[Theorem 3.7.]{Hofmann-Mitreas2017} are satisfied, once we confirm 1.,2.,3.: We choose a smooth decomposition of unity $b_Q \in C_c^\infty(\M)$ each supported within a coordinate patch of $\M$ and constantly one in a small Whitney cube $Q$. Then it suffices to show that for some $\alpha \in (0,1]$ the following holds for any $\eta \in C_c^\infty(B(x_0,r))$:
\begin{equation}\label{eq:Tetaest}
 \int_{B(x_0,r)} \int_0^r |T\eta(x,t)|^2 \frac{\mathrm dt}{t}\, \mathrm dx \aleq r^{n+2\alpha} [\eta]^2_{C^\alpha}.
\end{equation}
But observe that because of $\int_{\mathcal{M}} \kappa(x,y,t) \, \mathrm dy = 0$ by assumption (and $\kappa$ is integrable by the assumptions as well)
\[
\begin{split}
 |T\eta(x,t)| =& \left |\int_{\mathcal{M}} \kappa(x,y,t) \, (\eta(y)-\eta(x))\, \mathrm dy \right |\\
 \aleq& [\eta]_{C^\alpha}\, \int_{\mathcal{M}} \frac{t^{\nu}}{\brac{d(x,y)^2+t^2}^{\frac{n+\nu}{2}}}\, d(x,y)^\alpha\, \mathrm dy\\
 =& [\eta]_{C^\alpha}\, t^{-n+\alpha} \int_{\mathcal{M}} \frac{\brac{\frac{d(x,y)}{t}}^\alpha}{\brac{\brac{\frac{d(x,y)}{t}}^2+1}^{\frac{n+\nu}{2}}}\,  \mathrm dy\\
\end{split}
 \]
which holds for any $\alpha \in [0,1]$.
Now,
\[
\begin{split}
 \int_{\mathcal{M}} \frac{\brac{\frac{d(x,y)}{t}}^\alpha}{\brac{\brac{\frac{d(x,y)}{t}}^2+1}^{\frac{n+\nu}{2}}}\,  \mathrm dy	
 =&\sum_{k=1}^\infty \int_{B(x,2^k t) \backslash B(x,2^{k-1}t)} \frac{\brac{\frac{d(x,y)}{t}}^\alpha}{\brac{\brac{\frac{d(x,y)}{t}}^2+1}^{\frac{n+\nu}{2}}}\,  \mathrm dy\\
 &+\int_{B(x,t)} \frac{\brac{\frac{d(x,y)}{t}}^\alpha}{\brac{\brac{\frac{d(x,y)}{t}}^2+1}^{\frac{n+\nu}{2}}}\,  \mathrm dy.
\end{split}
 \]
Observe that 
\[
 \int_{B(x,t)} \frac{\brac{\frac{d(x,y)}{t}}^\alpha}{\brac{\brac{\frac{d(x,y)}{t}}^2+1}^{\frac{n+\nu}{2}}}\, \mathrm dy \aleq \int_{B(x,t)} 1\, \mathrm dy \lesssim t^n
\]
and
\[
 \begin{split}
 \int_{B(x,2^k t) \backslash B(x,2^{k-1}t)} \frac{\brac{\frac{d(x,y)}{t}}^\alpha}{\brac{\brac{\frac{d(x,y)}{t}}^2+1}^{\frac{n+\nu}{2}}}\, \mathrm dy 
 \aleq & (2^k t)^n  \frac{2^{\alpha k}}{\brac{2^{2k}+1}^{\frac{n+\nu}{2}}} \lesssim t^n 2^{k(\alpha -\nu)}.
\end{split}
 \]
We conclude that for $\alpha < \nu$,
\[
 |T\eta(x,t)| \aleq 
  [\eta]_{C^\alpha}\, t^{\alpha}.
 \]
This implies \eqref{eq:Tetaest}.
\end{proof}

The main point is that an admissible manifold allows for the $T(b)$-theorem to be applied to the $\sigma$-harmonic extension.

As a corollary we obtain a result which is essentially a square function estimate.
\begin{proposition}\label{co:squarebounded}
Let $\M$ be admissible and $0 < \sigma < 1$. Let $U$ be the $\sigma$-harmonic extension of $u \in C_c^\infty(\M)$, then
\begin{equation}\label{eq:sq:1}
 \int_{\mathcal{M} \times (0,\infty)} t\, |\partial_t U(x,t)|^2 \, \mathrm dx\,\mathrm dt \aleq \|u\|_{L^2(\mathcal{M})}^2,
\end{equation}
\begin{equation}\label{eq:sq:2}
 \int_{\mathcal{M} \times (0,\infty)} t \,|\nabla_{\kern -2pt x} U(x,t)|^2 \, \mathrm dx\, \mathrm dt \aleq \|u\|_{L^2(\mathcal{M})}^2.
\end{equation}
\end{proposition}
\begin{proof}
Let $p(x,y,s)$ be the heat kernel for $\mathcal{M}$. Then
\[
\begin{split}
 U(x,t) :=& \frac{1}{4^\sigma \Gamma(\sigma)} \int_0^\infty \int_{\mathcal{M}} p(x,y,s) \, u(y) \, \mathrm dy \, t^{2\sigma} e^{-\frac{t^2}{4s}} \, \frac{\mathrm ds}{s^{1+\sigma}}\\
 =& \frac{1}{4^\sigma \Gamma(\sigma)} \int_0^\infty \int_{\mathcal{M}} p(x,y,t^2 s) \, u(y) \, \mathrm dy \, e^{-\frac{1}{4s}} \, \frac{\mathrm ds}{s^{1+\sigma}}.
\end{split}
 \]
Regarding \eqref{eq:sq:1}, we have
\[
  t \, \partial_t U(x,t) = \frac{1}{4^\sigma \Gamma(\sigma)} \int_0^\infty \int_{\mathcal{M}} p(x,y,s) \, u(y) \, t \,\partial_t \brac{t^{2\sigma} e^{-\frac{t^2}{4s}}}\, \mathrm dy\, \frac{\mathrm ds}{s^{1+\sigma}}.
\]
Thus, for
\begin{align*}
\kappa(x,y,t) 	& := \frac{1}{4^\sigma \Gamma(\sigma)} \int_0^\infty p(x,y,s)\, t \, \partial_t \brac{t^{2\sigma} e^{-\frac{t^2}{4s}}} \, \frac{\mathrm ds}{s^{1+\sigma}}	\\
 				& = 2 \, \frac{t^{2\sigma}}{4^\sigma \Gamma(\sigma)} \int_0^\infty p(x,y,s) \, \brac{\sigma - \frac{t^2}{4s}} \, e^{-\frac{t^2}{4s}} \, \frac{\mathrm ds}{s^{1+\sigma}},
\end{align*}

and
\begin{equation}\label{eq:Toperator}
 Tu(x,t) := \int_{\mathcal{M}} \kappa(x,y,t) \, u(y)\, \mathrm dy,
\end{equation}
we have
	\begin{align*}
		\int_{\mathcal{M} \times (0,\infty)} t \, |\partial_t U(x,t)|^2 \, \mathrm dx\, \mathrm dt = \int_{\mathcal{M}} \int_{0}^{\infty} |Tu(x,t)|^2 \frac{\mathrm dt}{t} \, \mathrm dx.
	\end{align*}
Since for $u$ constant we know that $U(x,t)$ is also constant (see Appendix \ref{a:heatkernelextension}), one has $t \, \partial_t U \equiv 0$. We conclude that
\[
 \int_{\mathcal{M}} \kappa(x,y,t)\, \mathrm dy = 0 \quad \text{for all } x \in \mathcal{M},\, t > 0.
\]

It remains to establish the estimates \eqref{eq:TBkappaest} and \eqref{eq:TBkappalipest}, then the claim follows from \Cref{Tb}.

We estimate
\begin{align*}
	|\kappa(x,y,t)| \lesssim t^{2\sigma} \int_0^\infty p(x,y,s) \, \big|\sigma - \frac{t^2}{4s} \, \big| \, e^{-\frac{t^2}{4s}} \, \frac{\mathrm ds}{s^{1+\sigma}}.
\end{align*}
Since $\M$ was assumed to be admissible, we can use \eqref{c1} (after the transformation $s \longmapsto t^2 s$) to conclude \eqref{eq:TBkappaest}.
In order to show \eqref{eq:TBkappalipest}, we use the mean value theorem  to rewrite
	\begin{align*}
		\frac{|\kappa(x,y_1,t) - \kappa(x,y_2,t)|}{d(y_1,y_2)} 	& = |\nabla_{\kern-2pt y} \kappa (x,y,t)|	\\
																& \lesssim t^{2\sigma} \int_0^\infty |\nabla_{\kern-2pt y} p(x,y,s)| \, \big|\sigma - \frac{t^2}{4s} \, \big| \, e^{-\frac{t^2}{4s}} \, \frac{\mathrm ds}{s^{1+\sigma}}.
	\end{align*}
Again, since $\M$ was assumed to be admissible, we can use \eqref{c3} to deduce \eqref{eq:TBkappalipest}.

In order to derive the second estimate \eqref{eq:sq:2}, we argue similarly for a slightly different kernel $\kappa$. By the representation formula for $U$, we can rewrite
	\begin{align*}
		t \, \nabla_{\kern-2pt x} U(x,t) = \int_\M \nabla_{\kern-2pt x} p(x,y,s) \, u(y) \, \mathrm dy \, t^{2\sigma + 1} \, e^{-\frac{t^2}{4s}} \frac{\mathrm ds}{s^{1 + \sigma}},
	\end{align*}
so we define the operator	
	\begin{align*}
		Tu(x,t) = \int_\M \kappa(x,y,t) \, u(y) \, \mathrm dy
	\end{align*}
with the kernel
\[
\begin{split}\kappa(x,y,t) 	 :=\, & t^{2\sigma + 1} \int_0^\infty \nabla_{\kern-2pt x} p(x,y,s) \, e^{-\frac{t^2}{4s}} \frac{\mathrm ds}{s^{1 + \sigma}}	\\
						 =\, & t \int_0^\infty \nabla_{\kern-2pt x} p(x,y,t^2 s) \, e^{-\frac{1}{4s}} \frac{\mathrm ds}{s^{1 + \sigma}}.
\end{split}
						\]
Then again it holds $\int \kappa = 0$ and
	\begin{align*}
		\int_{\mathcal{M} \times (0,\infty)} t \, |\nabla_{\kern-2pt x} U(x,t)|^2 \, dx\, dt = \int_{\mathcal{M}} \int_{0}^{\infty} |Tu(x,t)|^2 \frac{\mathrm dt}{t} \, \mathrm dx.
	\end{align*}
So it suffices to show the estimates \eqref{eq:TBkappaest} and \eqref{eq:TBkappalipest}. This follows analogously to the first case by the admissibility of $\M$ and the mean value theorem.

\end{proof}

\begin{corollary}\label{cor:Tbprep}
Let $\M$ be admissible and $0 < \sigma < 1$. Let $U$ denote the $\sigma$-harmonic extension of $u$, then 
\[
  Tu(x,t) := t \, \nabla_{\kern-2pt (x,t)} U(x,t)
\]
satisfies the conditions of \Cref{Tb}. In particular, the following estimate holds for all functions $u \in C_c^\infty(\M)$:
\begin{align*}
	\int_\M \int_{(0,\infty)} t \, |\nabla_{\kern-2pt(x,t)}U(x,t)|^2 \, \mathrm dt \, \mathrm dx \lesssim \|u\|_{L^2(\M)}^2.
\end{align*}
\end{corollary}
\begin{proof}
 Follows immediately from Proposition 3.2.
\end{proof}
\section{BMO and Carleson measures: Proof of \texorpdfstring{\Cref{main}}{the main theorem}}\label{s:mainproof}
The proofs of \Cref{main} and \Cref{main2} consist in proving two directions. The easier one is \Cref{pr:direction1}, the more difficult one is  \Cref{pr:direction2}. 
\begin{proposition}\label{pr:direction1} Let $\M$ be admissible and $0 < \sigma < 1$. Let $U$ be the $\sigma$-harmonic extension of $u \in C_c^\infty(\M)$, then
\[
 \sup_{B} |B|^{-1} \int_{T(B)} |\nabla_{(x,t)} U(x,t)|^2\, t\, \mathrm dx\,\mathrm dt \aleq [u]_{BMO(\mathcal{M})}^2.
\]
\end{proposition}
\begin{proof}
We extend the argument from \cite[IV, \textsection 4.3, pp.158f]{Stein-1993}.

Fix any ball $B \subset \mathcal{M}$, and denote by $B^\ast$ the ball with twice the radius.
We decompose
\[
 \nabla_{(x,t)} U = \nabla_{(x,t)} U_1 + \nabla_{(x,t)} U_2 + \nabla_{(x,t)} U_3,
\]
where $U_i$ is the $\sigma$-harmonic extension of $u_i$, respectively, given as
\[
 u_1 := \chi_{B^\ast} (u-(u)_{B^\ast}),
\]
\[
 u_2 := (1-\chi_{B^\ast}) (u-(u)_{B^\ast}),
\]
\[
 u_3 := (u)_{B^\ast}.
\]
Observe that $U_3$ is constant and thus $\nabla U_3 = 0$.
Moreover,
\[
 |B|^{-1}\int_{T(B)} |\nabla_{(x,t)} U_1(x,t)|^2\, t\, \mathrm dx\,\mathrm dt  \leq |B|^{-1}\int_{\mathcal{M} \times (0,\infty)} |\nabla_{(x,t)} U_1(x,t)|^2\, t\, \mathrm dx\,\mathrm dt.
\]
In view of \Cref{cor:Tbprep},
\[
 |B|^{-1}\int_{T(B)} |\nabla_{(x,t)} U_1(x,t)|^2\, t\, \mathrm dx\,\mathrm dt \aleq |B|^{-1}\brac{\int_{B^\ast} |u-(u)_{B^\ast}|^2}
\]
Since $\mathcal{M}$ is supposed to be Ahlfors-regular,  $|B| \aeq |B^\ast|$ with uniform constants, one has by John-Nirenberg inequality \cite[(5.8)]{Ryansbook}, see also \cite{CW77},
\[
 |B^\ast |^{-1}\brac{\int_{B^\ast} |u-(u)_{B^\ast}|^2} \aleq [u]_{BMO}^2.
\]
This implies,
\[
 |B|^{-1}\int_{T(B)} |\nabla_{(x,t)} U_1(x,t)|^2\, t\, \mathrm dx\,\mathrm dt  \aleq [u]_{BMO}^2.
\]
It remains to estimate $U_2$.
As in the proof of \Cref{co:squarebounded},
\[
 |B|^{-1}\int_{T(B)} |\nabla_{(x,t)} U_2(x,t)|^2\, t\, \mathrm dx\,\mathrm dt  = |B|^{-1} \int_{T(B)} \left |\int_{\mathcal{M}} \kappa(x,y,t) u_2(y)\, \mathrm dy\right |^2\, \mathrm dx\, \frac{\mathrm dt}{t}.
\]
In view of \eqref{eq:TBkappaest}, for some given $\nu > 0$,
\[
|B|^{-1}\int_{T(B)} |\nabla_{(x,t)} U_2(x,t)|^2 \, t\, \mathrm dx\,\mathrm dt   \aleq |B|^{-1}\int_{T(B)}\int_{\mathcal{M}}\frac{t^{\nu}}{\brac{d(x,y)^2 + t^2}^{\frac{n+\nu}{2}}} |u_2(y)|\, \mathrm dy\, \mathrm dx\, \frac{\mathrm dt}{t}
\]
We denote with $B_{k}$, $k \in \N$, the ball concentric around $B$ but with radius $2^k$ times the radius of $B$.

Since $\supp u_2 \subset \mathcal{M} \backslash B^\ast$, we find that for any $x \in B$
\[
\begin{split}
 \int_{\mathcal{M}}\frac{t^\nu}{\brac{d(x,y)^2 + t^2}^{\frac{n+\nu}{2}}} |u_2(y)| \mathrm dy \aleq\sum_{k=0}^\infty \frac{t^\nu}{\brac{(2^{k}r)^2 + t^2}^{\frac{n+\nu}{2}}}  \int_{B_{k+1}\backslash B_{k}}|u(y)-(u)_{B}| \, \mathrm dy.
\end{split} 
\]
By triangle inequality and a telescoping sum,
\[
 \int_{B_{k+1}\backslash B_{k}}|u(y)-(u)_{B}| \, \mathrm dy \leq \int_{B_{k+1}}|u(y)-(u)_{B_{k+1}}| \, \mathrm dy + |B_{k+1}|\, \sum_{i=0}^{k}|(u)_{B_i} - (u)_{B_{i+1}}|.
\]
Using the doubling property of the measure and the definition of BMO,
we have 
\[
 \int_{B_{k+1}\backslash B_{k}}|u(y)-(u)_{B}| \mathrm dy \aleq |B_{k}|\, \brac{k+1}\, [u]_{BMO}.
\]
Consequently,
\[
\begin{split}
 \int_{\mathcal{M}}\frac{t^\nu}{\brac{d_{\mathcal{M}}(x,y)^2 + t^2}^{\frac{n+\nu}{2}}} |u_2(y)| \, \mathrm dy \lesssim \sum_{k=0}^\infty \frac{t^\nu}{\brac{(2^{k}r)^2 + t^2}^{\frac{n+\nu}{2}}}  |B_{k}|\brac{k+1}\, [u]_{BMO}.
\end{split} 
\]
This implies
\[
 |B|^{-1}\int_{T(B)} |\nabla_{(x,t)} U_2(x,t)|^2\, t\, \mathrm dx\,\mathrm dt \aleq A(r)\, [u]_{BMO},
\]
where
\[
 A(r) = \int_{0}^r \sum_{k=0}^\infty \frac{t^\nu}{\brac{(2^{k}r)^2 + t^2}^{\frac{n+\nu}{2}}}  (2^{k}r)^n\, \brac{k+1}\, \, \frac{\mathrm dt}{t}.
\]
By a substitution $t \mapsto rt$ we see that $A(r) = A(1)$, and 
\[
 A(1) \leq \brac{\int_0^1 t^{\nu-1}\, \mathrm dt}\cdot \brac{\sum_{k=0}^\infty 2^{-\nu k}\, \brac{k+1}} = C_\nu < \infty.
\]
This concludes the proof of \Cref{pr:direction1}. 
\end{proof}

\begin{proposition}\label{pr:direction2}
Let $\M$ be admissible and $0 < \sigma < 1$. Let $U$ be the $\sigma$-harmonic extension of $u \in C_c^\infty(\M)$, then
\[
 [u]_{BMO(\mathcal{M})}^2 \aleq \sup_{B} |B|^{-1} \int_{T(B)} |\nabla_{(x,t)} U(x,t)|^2\, t\, \mathrm dx\,\mathrm dt
\]
\end{proposition}

One technical ingredient in the proof of \Cref{pr:direction2} is the following observation (cf. \cite[(40) p.163]{Stein-1993}).

\begin{lemma}\label{la:extension}
Let $\Phi,U: \mathcal{M} \times (0,\infty) \to \R$ be the $\sigma-$harmonic extension of $\varphi $ and $u$, respectively. Then
\[
\begin{split}
\left | \int_{\mathcal{M}} u\, \varphi \right |\aleq& \int_{\mathcal{M}\times (0,\infty)} t |\partial_t \Phi(x,t)|\, |\partial_t U(x,t)|\, \mathrm dx\, \mathrm dt\\
 &+ \int_{\mathcal{M}\times (0,\infty)} t\, |\nabla_x \Phi(x,t)|\, |\nabla_x U(x,t)|\, \mathrm dx\, \mathrm dt
\end{split}
 \]
\end{lemma}
\begin{proof}
By integration by parts and the decay as $t \to \infty$ (see \Cref{a:heatkernelextension}), we have for every $x \in \mathcal{M}$:
\[
 u(x)\, \varphi(x) = \frac{1}{2\sigma} \int_{0}^\infty t^{2\sigma} \partial_{t} \brac{t^{1-2\sigma} \partial_t \brac{\Phi(x,t) U(x,t)}}\, \mathrm dt.
\]
Since $U$ is the $\sigma$-harmonic extension,
\[
\partial_{t} \brac{t^{1-2\sigma} \partial_t U(x,t)} = -t^{1-2\sigma} \lap_x U(x,t)
\]
and likewise for $\Phi$
%
after integration by parts
\begin{align*}
\int_\M u(x) \, \varphi(x) \ \mathrm dx = \, 	& \frac 1\sigma \int_{\M \times (0,\infty)} t \, \partial_t U(x,t) \, \partial_t \Phi(x,t) \ \mathrm dx \ \mathrm dt	\\
									& + \frac 1\sigma \, \int_0^\infty t \, \langle \nabla_{\kern-2pt x} U(\,\cdot\, , t), \nabla_{\kern-2pt x} \Phi(\, \cdot \, , t)\rangle_{L^2(\M)} \ \mathrm dt,
\end{align*}
from which the claim follows immediately.
\end{proof}

The following is proved in \cite[IV, \textsection 4.3, pp.162, Proposition]{Stein-1993}. It is stated there in $\R^n$, but the proof easily extends almost verbatim to Ahlfors regular spaces. 
\begin{lemma}\label{la:carlesonhoelder}
Let $\mathcal{M}$ be an Ahlfors-regular manifold and let $F,G \colon \M \times (0,\infty) \longrightarrow \R$ be measurable functions. Then,
\[
\begin{split}
 &\int_{\mathcal{M}\times (0,\infty)} t\, |F(x,t)| \, |G(x,t)| \, \mathrm dx\, \mathrm dt \\
 \aleq &\brac{\int_{\mathcal{M}} \brac{\int_{d(x,y) < t}|F(y,t)|^2\, \frac{\mathrm dy\, \mathrm dt}{t^{n-1}}}^{\frac{1}{2}}\, \mathrm dx} \brac{\sup_{B} |B|^{-1} \int_{T(B)} |G(x,t)|^2\, t\, \mathrm dx\,\mathrm dt}^{\frac{1}{2}}
\end{split}
 \]

\end{lemma}
\begin{proof}[Proof of \Cref{pr:direction2}]
Again, we essentially can follow Stein's book, namely \cite[IV, \textsection 4.3, pp.163f]{Stein-1993}.

From \Cref{la:extension} and \Cref{la:carlesonhoelder} for $F = \nabla_{\kern -2pt (x,t)} \Phi$ and $G = \nabla_{\kern -2pt (x,t)} U$ we obtain (using the notation \eqref{eq:carlesonbmo})
\[
 \left | \int_{\mathcal{M}} u\, \varphi\ \right | \aleq \brac{\int_{\mathcal{M}} \brac{\int_{d(x,y) < t}|\nabla \Phi(y,t)|^2\, \frac{\mathrm dy\, \mathrm dt}{t^{n-1}}}^{\frac{1}{2}}\, \mathrm dx}\, [u]_{\tilde{BMO}}.
\]
We can conclude once we show that
\begin{equation}\label{eq:hardymaxest}
\brac{\int_{\mathcal{M}} \brac{\int_{d(x,y) < t}|\nabla \Phi(y,t)|^2\, \frac{\mathrm dy\, \mathrm dt}{t^{n-1}}}^{\frac{1}{2}}\, \mathrm dx} \aleq \|\varphi\|_{\mathcal{H}^1},
\end{equation}
where $\mathcal{H}^1$ is the Hardy space. Indeed, the claim then follows in view of the duality of Hardy spaces and BMO \cite[(7.154)]{Ryansbook}.

To obtain \eqref{eq:hardymaxest} we use an extrapolation result \cite[Theorem 6.18.]{Hofmann-Mitreas2017}, which essentially states that a suitable operator, if it is bounded from $L^2$ to $L^2$, can be extended to a bounded operator from the Hardy space into $L^1$. 
To apply this result, first observe that from Fubini we have
\begin{equation} \label{eq:cone}
 \begin{split}
  &\left \|x \mapsto \brac{\int_{d(x,y) < t}|\nabla \Phi(y,t)|^2\, \frac{\mathrm dy\, \mathrm dt}{t^{n-1}}}^{\frac{1}{2}} \right \|_{L^2(\M)}^2\\
  =& \int_0^\infty \int_{\M} \int_{\M} t^{1-n} \chi_{d(x,y) < t} |\nabla \Phi(y,t)|^2\, \mathrm dx\, \mathrm dy\, \mathrm dt\\
  \aleq & \int_0^\infty \int_{\M}  t  |\nabla \Phi(y,t)|^2\, \mathrm dy\, \mathrm dt.
 \end{split}
\end{equation}
From  \Cref{cor:Tbprep} we conclude that the operator $\varphi \mapsto T\varphi := t \nabla \Phi$ as in \eqref{eq:Toperator} satisfies the conditions imposed on the operator denoted by $\theta_E$ in \cite[Theorem 6.18.]{Hofmann-Mitreas2017}. That theorem implies \eqref{eq:hardymaxest}.
\end{proof}

\section{Jacobian Estimate: Proof of Theorem~\ref{th:jac}}\label{s:thjac}
By the trace space characterizations of BMO obtained in Theorem~\ref{main} and \Cref{main2} we can follow the strategy in \cite{LS} to prove the Coifman--Lions--Meyer--Semmes estimate on manifolds.

\begin{proof}[Proof of Theorem~\ref{th:jac}]
For some $\frac 12 < \sigma < 1$, let $F^\ell$ be the $\sigma$-harmonic extension to $\M \times (0,\infty)$ for $f^\ell$ for all $\ell = 1, \ldots, n$, and let $\Phi$ be the $\sigma$-harmonic extension of $\varphi$. Then by Stokes theorem we have
	\begin{align*}
		\left |\int_\M \mathrm df^1\wedge \ldots \wedge \mathrm df^n\ \varphi \right |= \left |\int_{\M \times (0,\infty)}  \mathrm d_{(x,t)} F^1 \wedge \ldots \wedge \mathrm d_{(x,t)} F^n \wedge \mathrm d_{(x,t)} \Phi \right | .
	\end{align*}
We claim that 
\begin{equation}\label{eq:jac:claim}
 \begin{split}&\left |\int_{\M \times (0,\infty)}  \mathrm d_{(x,t)} F^1 \wedge \ldots \wedge  \mathrm d_{(x,t)} F^n \wedge \mathrm d_{(x,t)} \Phi \right | \\
\aleq& \int_{\M \times (0,\infty)} t\, |\nabla_{(x,t)} F|^{n-1}\, |\nabla_x \nabla_{(x,t)} F|\, |\nabla_{(x,t)}\Phi|.\\
\end{split}
\end{equation}
Assume we have \eqref{eq:jac:claim} (which will be proven below). Set
\[
 M_f(x) := \sup_{d(x,y) < t} |\nabla_{(x,t)} F(y,t)|.
\]
Then, by Lemma~\ref{la:carlesonhoelder}, H\"older's inequality, and the BMO-characterization, Theorem~\ref{main}
\[
\begin{split}
 &\int_{\mathcal{M} \times (0,\infty)} t\,|\nabla_{(x,t)} F|^{n-1}\, |\nabla_x \nabla_{(x,t)} F|\, |\nabla_{(x,t)}\Phi| \\
 \aleq & \brac{\int_{\mathcal{M}} \brac{\int_{d(x,y) < t}\brac{|\nabla_{(x,t)} F|^{n-1}\, |\nabla_x \nabla_{(x,t)} F(y,t)|}^2\, \frac{\mathrm dy\, \mathrm dt}{t^{n-1}}}^{\frac{1}{2}}\, \mathrm dx} [\varphi]_{BMO}\\
 \aleq & \int_{\mathcal{M}} |M_f|^{n-1} \brac{\int_{d(x,y) < t}\brac{\, |\nabla_x \nabla_{(x,t)} F(y,t)|}^2\, \frac{\mathrm dy\, \mathrm dt}{t^{n-1}}}^{\frac{1}{2}}\, \mathrm dx\ [\varphi]_{BMO}\\
 \aleq & \|M_f\|_{L^n(\M)}^{n-1} \left \| \brac{\int_{d(\cdot,y) < t}\brac{\, |\nabla_x \nabla_{(x,t)} F(y,t)|}^2\, \frac{\mathrm dy\, \mathrm dt}{t^{n-1}}}^{\frac{1}{2}} \right \|_{L^n(\M)}  [\varphi]_{BMO}.
\end{split}
 \]
By Lemma~\ref{la:maximalfct}, and the boundedness of the maximal function from $L^n(\M)$ to $L^n(\M)$ (which holds true on every space with doubling measure, \cite{Stein-1993})
\[
 \|{M}_f\|_{L^n(\M)} \aleq \|\nabla f\|_{L^n(\M)}.
\]
Moreover, similar as in \eqref{eq:cone} together with the extrapolation result \cite[Theorem 6.18]{Hofmann-Mitreas2017}, we obtain
\[
 \left \| \brac{\int_{d(\cdot,y) < t}\brac{\, |\nabla_x \nabla_{(x,t)} F(y,t)|}^2\, \frac{\mathrm dy\, \mathrm dt}{t^{n-1}}}^{\frac{1}{2}} \right \|_{L^n(\M)} 
 \aleq \|\nabla f\|_{L^n(\M)}.
\]
Here we also have used that $\nabla_x F$ is the $\sigma$-harmonic extension of $\nabla f$, which follows from the uniqueness of the $\sigma$-harmonic extension.

This concludes the proof of Theorem~\ref{th:jac}, up to proving \eqref{eq:jac:claim}.
\end{proof}

\begin{proof}[Proof of \eqref{eq:jac:claim}]
 
	Using one integration by parts in $t$-direction, together with the decay of the $\sigma$-harmonic extension for $t \to \infty$, we obtain
	\[
	\begin{split}
		&\left |\int_{\M \times (0,\infty)} \mathrm d_{(x,t)} F^1 \wedge \ldots \wedge \mathrm d_{(x,t)} F^n \wedge \mathrm d_{(x,t)} \Phi\right | \\
		=& \left |\int_{\M \times (0,\infty)} t^{2\sigma} \, \partial_t [ t^{1 - 2\sigma} \mathrm d_{(x,t)} F^1 \wedge \ldots \wedge \mathrm d_{(x,t)} F^n \wedge \mathrm d_{(x,t)} \Phi] \right |.
		\end{split}
	\]
By orthogonility of the coordinates in $\M \times (0,\infty)$ we have \[\mathrm d_{(x,t)} F(x,t)= \mathrm d_{\M} F(x,t) + \partial_t F(x,t) \, \mathrm dt.\]
Since $\mathrm dt \wedge \mathrm dt = 0$, we have two cases to estimate:
\begin{equation}\label{eq:jac:1}
 \left |\int_{\M \times (0,\infty)} t^{2\sigma} \, \partial_t [ t^{1 - 2\sigma} \mathrm d_{\mathcal{M}} F^1\wedge \ldots \wedge \partial_t F^\ell\, \mathrm dt\wedge \ldots \wedge \mathrm d_{\mathcal{M}} F^n \wedge \mathrm d_{\mathcal{M}} \Phi] \right |
\end{equation}
for $\ell = 1,\ldots,n$ and
\begin{equation}\label{eq:jac:2}
 \left |\int_{\M \times (0,\infty)} t^{2\sigma} \, \partial_t [ t^{1 - 2\sigma} \mathrm d_{\mathcal{M}} F^1\wedge \ldots \wedge \mathrm d_{\mathcal{M}} F^n \wedge \partial_t \Phi\, \mathrm dt] \right |.
\end{equation}
Regarding \eqref{eq:jac:1}, 
\[
\begin{split}
 & \left |\int_{\M \times (0,\infty)} t^{2\sigma} \, \partial_t [ t^{1 - 2\sigma} \mathrm d_{\mathcal{M}} F^1\wedge \ldots \wedge \partial_t F^\ell\, \mathrm dt\wedge \ldots \wedge \mathrm d_{\mathcal{M}} F^n \wedge \mathrm d_{\mathcal{M}} \Phi] \right |\\
\leq & \left |\int_{\M \times (0,\infty)} t\, \mathrm d_{\mathcal{M}} F^1\wedge \ldots \wedge \partial_t F^\ell\, \mathrm dt\wedge \ldots \wedge \mathrm d_{\mathcal{M}} F^n \wedge \mathrm d_{\mathcal{M}} \partial_t\Phi \right |\\
&+ \int_{\M \times (0,\infty)} t^{2\sigma} |\nabla_{(x,t)} F|^{n-1}\, |\partial_t(t^{1-2\sigma})\partial_t F| |\nabla_x \Phi|\\
&+\int_{\M \times (0,\infty)} t |\partial_t \nabla_x F| |\nabla_{(x,t)} F|^{n-1} |\nabla_x	 \Phi|.
\end{split}
\]
The last term is already in the form of \eqref{eq:jac:claim}. The second-to-last term is as well, if we use that $\partial_t(t^{1-2\sigma}\partial_t F) = t^{1-2\sigma}\lap_{\M} F$.
For the first term we use an integration by parts (observe that $\M$ has no boundary)
\[
\begin{split}
 &\left |\int_{\M \times (0,\infty)} t\, \mathrm d_{\mathcal{M}} F^1\wedge \ldots \wedge \partial_t F^\ell\, \mathrm dt\wedge \ldots \wedge \mathrm d_{\mathcal{M}} F^n \wedge \mathrm d_{\mathcal{M}} \partial_t\Phi \right |\\
 =&\left |\int_{\M \times (0,\infty)} t\, \mathrm d_{\mathcal{M}}\brac{\mathrm d_{\mathcal{M}} F^1\wedge \ldots \wedge \partial_t F^\ell\, \mathrm dt\wedge \ldots \wedge \mathrm d_{\mathcal{M}} F^n} \partial_t\Phi \right |\\
 \end{split}
\]
Now it follows from Leibniz rule that this is of the form of \eqref{eq:jac:claim}, so \eqref{eq:jac:1} can be estimated as claimed.

Regarding the estimate for the term \eqref{eq:jac:2} we argue similarly.
\[
\begin{split}
  &\left |\int_{\M \times (0,\infty)} t^{2\sigma} \, \partial_t [ t^{1 - 2\sigma} d_{\mathcal{M}} F^1 \wedge \ldots \wedge \mathrm d_{\mathcal{M}} F^n \wedge \partial_t \Phi\, \mathrm dt] \right |\\
  \aleq&\left |\int_{\M \times (0,\infty)} t^{2\sigma} \, d_{\mathcal{M}} F^1 \wedge \ldots \wedge \mathrm d_{\mathcal{M}} F^n \wedge \partial_t [ t^{1 - 2\sigma} \partial_t \Phi]\, \mathrm dt \right |\\
  &+\int_{\M \times (0,\infty)} t |\nabla_{(x,t)} \nabla_x F|\, |\nabla_x F|^{n-1} |\partial_t \Phi|\\
\end{split}
  \]
The second term is already of the form \eqref{eq:jac:claim}. For the first term we use again that $\partial_t(t^{1-2\sigma} \partial_t \Phi) = t^{1-2\sigma} \lap_{\M} \Phi$. An integration by parts along $\M$ (as we did for the estimate of \eqref{eq:jac:1} above) results again in an estimate of the form \eqref{eq:jac:claim}.
\eqref{eq:jac:2} is now estimated, and we conclude that \eqref{eq:jac:claim} holds true.
\end{proof}

\appendix
\section{Computations for the \texorpdfstring{$\sigma$}{sigma}-harmonic extension}\label{a:heatkernelextension}
We recall our definition for the extension $U$, motivated by \cite{Stinga-Torrea-2010}. 
Given $u \in C_c^\infty(\mathcal{M})$ and $0 < \sigma < 1$, the $\sigma$-harmonic extension $U: \mathcal{M} \times [0,\infty) \to \R$ is formally given by
\[
 U(x,t) := \frac{1}{4^\sigma \Gamma(\sigma)} t^{2\sigma} \int_0^\infty e^{s\lap_{\mathcal{M}}} u(x)\, e^{-\frac{t^2}{4s}}\, \frac{\mathrm ds}{s^{1+\sigma}}. 
\]
More explicitly, one has
\[
\begin{split}
 U(x,t) :=& \frac{1}{4^\sigma \Gamma(\sigma)} \int_0^\infty \int_{\mathcal{M}} p(x,y,s) \, u(y) \, \mathrm dy\, t^{2\sigma} e^{-\frac{t^2}{4s}}\, \frac{\mathrm ds}{s^{1+\sigma}},	\\
 			= & \frac{1}{\Gamma(\sigma)} \int_0^\infty \int_\M p(x,y,\tfrac{t^2}{4s}) \, u(y) \, \mathrm dy \, e^{-s} \, s^{\sigma - 1} \, \mathrm ds,
\end{split}
 \]
where $p(x,y,s)$ is the heat kernel for $\mathcal{M}$. Furthermore, $U$ is the smooth (in the interior) solution of
\begin{equation}\label{eq:Upde}
 \begin{cases}
 \lap_{\mathcal{M}} U  + \frac{1-2\sigma}{t} \partial _t U + \partial_{tt} U = 0 \quad& \mbox{in $\mathcal{M}\times (0,\infty)$}\\
 U(x,0) = u(x) & \mbox{in $\mathcal{M}$}\\
\lim_{|(x,t)| \to \infty} U(x,t) = 0. 
 \end{cases}
\end{equation}
The most important property for us is that constants are extended by constants: our manifold being assumed to be Ahlfors regular, it is stochastically complete, i.e
	\begin{align*}
		\int_\M p(x,y,s) \, \mathrm dy = 1	&& \text{for all } x \in \M, \ s > 0,
	\end{align*}
see \cite[Theorem 1]{G86}. Now using the representation formula, we deduce that for constant $u \colon \M \longrightarrow \mathbb R$
	\begin{align*}
		U(x,t)	& = \frac{u}{\Gamma(\sigma)} \int_0^\infty \int_\M p(x,y,\tfrac{t^2}{4s}) \, \mathrm dy \, e^{-s} s^{\sigma - 1} \, \mathrm ds = u. \\
	\end{align*}
One can also check, that the representation formula solves the PDE \eqref{eq:Upde}. Using the second line in the representation formula and that $p$ solves the heat equation \eqref{eq:heatkerneleq} 
 one sees that
	\begin{align*}
		\partial_t U(x,t)	& = \frac{1}{\Gamma(\sigma)} \int_0^\infty \int_\M (\partial_t p)(x,y,\tfrac{t^2}{4s}) \, u(y) \, \mathrm dy \, e^{-s} \, \Big[\frac t2 \, s^{\sigma - 2}\Big] \, \mathrm ds	\\
		\partial_t^2 U(x,t)	& = \frac{1}{\Gamma(\sigma)} \int_0^\infty \int_\M (\partial_t p)(x,y,\tfrac{t^2}{4s}) \, u(y) \, \mathrm dy \, e^{-s} \, \Big[-\frac 12 \, s^{\sigma - 2} - s^{\sigma - 1} - (1 - \sigma)s^{\sigma - 2} \Big] \, \mathrm ds	\\
		\Delta_x U(x,t) 	& =  \frac{1}{\Gamma(\sigma)} \int_0^\infty \int_\M (\partial_t p)(x,y,\tfrac{t^2}{4s}) \, u(y) \, \mathrm dy \, e^{-s} \, \Big[s^{\sigma - 1} \Big] \, \mathrm ds.
	\end{align*}
Comparing the square brackets shows that $U$ indeed solves the PDE \eqref{eq:Upde}.

Moreover, since the heat kernel is an approximation of the identity, we see that $U$ has also the correct boundary data. Even more, as $t \to \infty$, we have $U \to 0$. This follows from the admissibility of $\M$, i.e. given $u \in C_c^\infty
(\M)$, we can compute
	\begin{align*}
		 |U(x,t)|	& \lesssim t^{2\sigma} \int_{\mathcal{M}} |u(y)| \int_0^\infty p(x,y,s) \, e^{-\frac{t^2}{4s}} \, \frac{\mathrm ds}{s^{1+\sigma}} \, \mathrm dy \\
		 			& = \int_{\mathcal{M}} |u(y)| \int_0^\infty p(x,y,t^2 s) \, e^{-\frac{1}{4s}} \, \frac{\mathrm ds}{s^{1+\sigma}} \, \mathrm dy	\\
		 			& \lesssim \int_\M |u(y)| \, \frac{t^\nu}{(d(x,y)^2 + t^2)^{\frac{n + \nu}{2}}} \, \mathrm dy	\\
		 			& \leq \frac{1}{t^n} \, \int_\M |u(y)| \, \mathrm dy.
	\end{align*}
With an analogue computation, we see that
	\[
		|\partial_t U(x,t)| \lesssim t^{-(n + 1)} \, \int_\M |u(y)| \,  \mathrm dy.
	\]

We also record the following estimate by maximal functions of the $\sigma$-harmonic extension.
\begin{lemma}\label{la:maximalfct}
Let $\M$ be as in Theorem~\ref{main} or Theorem~\ref{main2}, and let $U$ denote the $\sigma$-harmonic extension of $u$, for $\frac 12 < \sigma < 1$. Then,
\[
 \sup_{y,t:\, d(x,y) < t} |\nabla_{(y,t)} U(y,t)| \aleq\, M|\nabla u|(x)
\]
where $M$ denotes the Hardy-Littlewood maximal function	
	\begin{align*}
		Mf(x) = \sup_{r > 0} \frac{1}{|B(x,r)|} \int_{B(x,r)} |f(y)| \, \mathrm dy.
	\end{align*}
\end{lemma}
\begin{remark}The estimate in Lemma~\ref{la:maximalfct} is false for $\sigma \leq \frac{1}{2}$ even in Euclidean space. E.g. for $\sigma = \frac{1}{2}$ we have
\[
 |\laps{1} u(x)| = \lim_{t \to 0^+} |\partial_t U(x,t)| \not \aleq M|\nabla_x u|(x).
\]
Indeed, the latter inequality has to be false because otherwise Lipschitz functions (i.e. functions with a finite right-hand side in the estimate above) have half-Laplacian bounded, which is false in general.\end{remark}
\begin{proof}[Proof of Lemma~\ref{la:maximalfct}]
	Let $p(x,y,s)$ be the heat kernel of $\M$, i.e. it solves \eqref{eq:heatkerneleq}. Moreover, let $t > 0$, $y \in \M$ and $x \in B(y,t)$. We will estimate $\nabla_y U$ and $\partial_t U$ seperately. Using the representation formula we see that 
		\begin{align*}
			\partial_t U(y,t)	& = \frac{1}{\Gamma(\sigma)} \int_0^\infty \int_\M (\partial_t p)(y,z,\tfrac{t^2}{4s}) \, u(z) \, \mathrm dz \, e^{-s} \, \Big[\frac t2 \, s^{\sigma - 2}\Big] \, \mathrm ds	\\
								& = \frac{1}{\Gamma(\sigma)} \int_0^\infty \int_\M \Delta_z p(y,z,\tfrac{t^2}{4s}) \, u(z) \, \mathrm dz \, e^{-s} \, \Big[\frac t2 \, s^{\sigma - 2}\Big] \, \mathrm ds	\\
								& = -\frac{1}{\Gamma(\sigma)} \int_0^\infty \int_\M \langle \nabla_z p(y,z,\tfrac{t^2}{4s}), \nabla u(z)\rangle \, \mathrm dz \, e^{-s} \, \Big[\frac t2 \, s^{\sigma - 2}\Big] \, \mathrm ds.
		\end{align*}
	Taking the absolute value and using Fubini we arrive at
		\begin{align*}
			|\partial_t U(y,t)| \lesssim \int_\M |\nabla u(z)| \int_0^\infty t \, |\nabla_z p(y,z,\tfrac{t^2}{4s})| \, e^{-s} \, s^{\sigma - 2} \mathrm ds \, \mathrm dz.
		\end{align*}
	Using the gradient estimate \eqref{eq:heat2} of the heat kernel, we deduce
		\begin{align*}
			\int_0^\infty t \, |\nabla_y p(y,z,\tfrac{t^2}{4s})| \, e^{-s} \, s^{\sigma - 2} \mathrm ds	& \lesssim \frac{1}{t^n} \int_0^\infty s^{\frac {n + 1}2 + \sigma - 2} \, e^{(-c \frac{d(y,z)^2}{t^2} + 1)s} \mathrm ds	\\
																										& \lesssim \frac{1}{t^n} \, \frac{1}{(c \frac{d(y,z)^2}{t^2} + 1)^{\frac n2 + \sigma - \frac 12}}.
		\end{align*}
	Let $B = B(x,2t)$ and $B_k := B(x,2^k t) \setminus B(x,2^{k - 1} t)$ for $k \geq 2$. Then we can estimate
		\begin{align} \label{eq:annuli_arg}
			|\partial_t U(y,t)|	&  \lesssim \frac{1}{t^n}\int_B  \, \frac{|\nabla u(z)| \, \mathrm dz}{(c \frac{d(y,z)^2}{t^2} + 1)^{\frac n2 + \sigma - 1}} + \sum_{k = 2}^\infty \frac{1}{t^n}\int_{B_k} \frac{|\nabla u(z)| \, \mathrm dz}{(c \frac{d(y,z)^2}{t^2} + 1)^{\frac n2 + \sigma - \frac 12}}
		\end{align}
	By the choice of the annuli, we can estimate for $z \in B_k$ and $k \geq 2$ 
		\begin{align*}
			d(y,z) \geq |d(x,z) - d(x,y)| \geq (2^{k - 1} - 1) t.
		\end{align*}
	So it follows that
		\[
		\begin{split}
			\sup_{d(x,y) < t} |\partial_t U(y,t)| & \lesssim \frac{1}{t^n} \int_B |\nabla u(z)| \, \mathrm dz + \sum_{k = 2}^\infty \frac{1}{2^{2k( \sigma - \frac 12)}} \, \frac{1}{(2^k t)^n} \int_{B_k} |\nabla u(z)| \, \mathrm dz\\
			\lesssim &M(|\nabla u|)(x),
		\end{split}
		\]
	since $\sigma > \frac 12$ and $\M$ is Ahlfors regular.
	
The estimate of $|\nabla_y U|$ is simpler (and works for any $0 < \sigma < 1)$. We simply observe that
		\begin{align*}
			(\nabla_y U)(y,t) 	& = \frac{1}{\Gamma(\sigma)} \int_0^\infty \int_\M \nabla_y p(y,z,\tfrac{t^2}{4s}) \, u(z) \, \mathrm dz \, e^{-s} \, s^{\sigma - 1} \, \mathrm ds	\\
								& = \frac{1}{\Gamma(\sigma)} \int_0^\infty \int_\M \nabla_z p(y,z,\tfrac{t^2}{4s}) \, u(z) \, \mathrm dz \, e^{-s} \, s^{\sigma - 1} \, \mathrm ds	\\
								& = -\frac{1}{\Gamma(\sigma)} \int_0^\infty \int_\M p(y,z,\tfrac{t^2}{4s}) \, \nabla u(z) \, \mathrm dz \, e^{-s} \, s^{\sigma - 1} \, \mathrm ds	\\
								& = -\frac{1}{\Gamma(\sigma)} \int_\M \nabla u(z) \int_0^\infty p(y,z,\tfrac{t^2}{4s}) \, e^{-s} \, s^{\sigma - 1} \ \mathrm ds \, \mathrm dz.
		\end{align*}
	Taking the absolute value we arrive at
		\begin{align*}
			|\nabla_y U(y,t)| 	& \lesssim \int_\M |\nabla u(z)| \int_0^\infty p(y,z,\tfrac{t^2}{4s}) \, e^{-s} \, s^{\sigma - 1} \, \mathrm ds \, \mathrm dz.
		\end{align*}
	Using \eqref{eq:heat1} we see that
		\begin{align*}
			\int_0^\infty p(y,z,\tfrac{t^2}{4s}) \, e^{-s} \, s^{\sigma - 1} \, \mathrm ds \lesssim \frac{1}{t^n} \, \frac{1}{(\frac{d(y,z)^2}{t^2} + 1)^{\frac n2 + \sigma}}.
		\end{align*}
	In the same fashion as before we can obtain (for any $0 < \sigma < 1$) the estimate
		\[
			\sup_{d(y,x) < t}|\nabla_y U(y,t)| \lesssim M(|\nabla u|)(x).
		\]
This concludes the proof.
\end{proof}

\bibliographystyle{abbrv}%
\bibliography{bib}%

\end{document}